\documentclass[11pt]{article}
\pdfoutput=5
\date{}
\makeatletter
\@addtoreset{equation}{section}
\makeatother
\usepackage{amssymb,amsmath}
\usepackage{cite}
\usepackage{hyperref}
\begin{document}

\title{\bf Normalized bound state solutions of fractional Schr\"{o}dinger equations with general potential\footnote{Supported by National Natural Science Foundation of China(No. 12271443) and Natural Science Foundation of Chongqing, China(cstc2020jcyjjqX0029)}}
\author{{Xin Bao,\ Ying Lv,\ Zeng-Qi Ou \footnote{Corresponding author.
E-mail address: bxinss@163.com(X. Bao), ly0904@swu.edu.cn(Y. Lv), ouzengqi@swu.edu.cn(Z.-Q. Ou).}}\\
{\small \emph{ School  of  Mathematics  and  Statistics, Southwest University,  Chongqing {\rm400715},}}\\
{\small \emph{People's Republic of China}}}
\maketitle
\baselineskip 17pt

{\bf Abstract}: In this paper, we study a class of fractional Schr\"{o}dinger equation
\begin{equation} \label{eq0}
\left\{
\begin{aligned}
&(-\Delta)^{s}u=\lambda u+a(x)|u|^{p-2}u,\\
&\int_{\mathbb{R}^{N}}|u|^{2}dx=c^{2},\ u\in H^{s}(\mathbb{R}^{N}),
\end{aligned}
\right.
\end{equation}
where $N>2s$, $s\in(0,1)$ and $p\in(2,2+4s/N), c>0$. $a(x)\in C(\mathbb{R}^{N},\mathbb{R})$ is a positive potential function. By using Fixed Point Theorem of Brouwer, barycenter function and variational method, we obtain the existence of normalized bound solutions for problem \eqref{eq0}.

\par

\section{Introduction and main results}
\ \ \ \ \ In this paper, we discuss the existence of positive solutions for a class of fractional Schr\"{o}dinger equation
\begin{equation} \label{eq1}
\left\{
\begin{aligned}
&(-\Delta)^{s}u=\lambda u+a(x)|u|^{p-2}u,\\
&\int_{\mathbb{R}^{N}}|u|^{2}dx=c^{2},\ u\in H^{s}(\mathbb{R}^{N}),
\end{aligned}
\right.
\end{equation}
where $N>2s$, $s\in(0,1)$ and $p\in(2,2+4s/N), c>0$. $a(x)\in C(\mathbb{R}^{N},\mathbb{R})$ satisfies the following assumption:

$(A_{1})$\ $0<\inf\limits_{x\in \mathbb{R}^N}a(x)\leq a(x)\leq\lim \limits_{|x|\rightarrow \infty} a(x)=a_\infty, a(x)\not\equiv a_\infty$ in $\mathbb{R}^{N}$.

 The fractional Laplacian $(-\Delta)^{s}$ is the nonlocal operator and defined by
$$
(-\Delta)^{s} u(x)=C_{N, s} \ P. V.\int_{\mathbb{R}^{N}} \frac{u(x)-u(y)}{|x-y|^{N+2 s}} dy,
$$
where $u:\mathbb{R}^{N}\rightarrow \mathbb{R}$ is sufficiently smooth, and $P.V.$ is an abbreviation for the Cauchy principal value and $C_{N, s}$ stands for a positive constant depending on $N$ and $s$(see \cite{4}).

For $s\in(0,1)$ and $N>2s$, the fractional Sobolev space of order $s$ on $\mathbb{R}^{N}$ is defined by
$$
H^s(\mathbb{R}^{N})=\left\{u \in L^2(\mathbb{R}^{N}): \int_{\mathbb{R}^{N}} \int_{\mathbb{R}^{N}} \frac{|u(x)-u(y)|^{2}}{|x-y|^{N+2 s}}dxdy<\infty\right\},
$$
endowed with the norm
$$
\|u\|=\left(\int_{\mathbb{R}^{N}}|(-\Delta)^{s / 2} u|^{2}dx+\int_{\mathbb{R}^{N}} u^{2}dx\right)^{1/2},
$$
where
$$
\int_{\mathbb{R}^{N}}|(-\Delta)^{\frac{s}{2}} u|^{2}d x=\int_{\mathbb{R}^{N}}\int_{\mathbb{R}^{N}} \frac{|u(x)-u(y)|^{2}}{|x-y|^{N+2 s}}dxdy.
$$
From \cite{4}, the embedding $H^s(\mathbb{R}^{N})\hookrightarrow L^q(\mathbb{R}^{N})$ is continuous for any $q\in [2,2^*_s]$, where $2^*_s=\frac{2N}{N-2s}$ is the fractional critical Sobolev exponent. Moreover, there is a positive constant $C=C(q,N,s)$ such that
\begin{equation} \label{eqA2}
\|u\|_{q}\leq C \|u\|
\end{equation}
for any $u\in H^s(\mathbb{R}^{N})$, where $\|\cdot\|_q$ denotes the norm of the Lebesgue space $L^q(\mathbb{R}^{N})$.

Now, we recall the Gagliardo-Nirenberg-Sobolev inequality(see \cite{8}): for any $p \in(2,2^{*}_s)$, there exists a sharp constant $C_{N, p,s}>0$ such that
\begin{equation}\label{eq6}
\|u\|^{p}_{p}\leq C_{N,p,s}\left(\int_{\mathbb{R}^{N}}|(-\Delta)^{s/2} u|^{2}dx\right)^{\frac{N(p-2)}{4s}}\|u\|_{2}^{p-\frac{N(p-2)}{2s}}
\end{equation}
for any $u \in H^{s}(\mathbb{R}^{N})$, where the best constant $C_{N, p,s}$ can be achieved.

In the past years, many researchers have extensively studied the existence of nontrivial solutions for the Sch\"{o}dinger equation with a prescribed $L^2$-norm:
\begin{equation} \label{eq2}
\left\{
\begin{aligned}
&-\Delta u+\lambda u=g(u),\quad x\in  \mathbb{R}^{N},\\
&\int_{\mathbb{R}^{N}}|u|^{2}dx=c^{2},
\end{aligned}
\right.
\end{equation}
 and such solution is known as a normalized solution and $\lambda$ appears as a Lagrange multiplier (see \cite{f}, \cite{Bartsch}, \cite{d}, \cite{CT}, \cite{DZ}, \cite{10}, \cite{JL}, \cite{MR}, \cite{e}, \cite{18} and the references therein). In the classical paper \cite{10}, when $g(u)$ is nonhomogeneous $L^{2}$-supercritical nonlinearity, by developing a mountain-pass argument on $S_c=\{u\in H^1(\mathbb{R}^{N}): \|u\|_2=c\}$, Jeanjean obtained the existence of a ground state normalized solution for problem \eqref{eq2}. When $N=2$ and $f$ satisfies an exponential critical growth in problem \eqref{eq2}, Chang, Liu and Yan in \cite{d} proved the existence of normalized ground state solution for any $c>0$ by using a constrained minimization method and the Trudinger-Moser inequality in $\mathbb{R}^2$.  In a mass subcritical case but under general assumptions on the nonlinearity $g$, Jeanjean and Lu in\cite{JL} obtained multiple nonradial solutions for problem \eqref{eq2} when $N=4$ or $N\geq 6$, and all these solutions are sign-changing.

For the following Sch\"{o}dinger equation with potential
\begin{equation} \label{eq2A}
\left\{
\begin{aligned}
&-\Delta u+V(x)u+\lambda u=g(u),\quad x\in  \mathbb{R}^{N},\\
&\int_{\mathbb{R}^{N}}|u|^{2}dx=c^{2},
\end{aligned}
\right.
\end{equation}
if $V\leq 0$, the existence of ground state normalized solution for problem \eqref{eq2A} is obtained in \cite{DZ}, where $g$ satisfies some Ambrosetti-Rabinowitz type conditions, and in \cite{MR} for the case  $g(u)=|p|^{p-2}u\ (p\in(2,2^*))$. But for the case $V\geq 0$, the problem becomes more delicate and difficult. When $V\geq 0$ satisfies $V(x)\to 0$ as $|x|\to \infty$ and $g(u)=|p|^{p-2}u\ (p\in(2+\frac{4}{N}, 2^*))$, by constructing a suitable linking geometry, Bartsch et al. \cite{Bartsch} succeeded to obtain the existence of bound state solution with high Morse index for problem \eqref{eq2A}.

In\cite{CT}, Chen and Tang proved the existence of normalized ground state solutions for the
following Schr\"{o}dinger equation
\begin{equation*}
\left\{
\begin{aligned}
&-\Delta u=\lambda u+a(x)g(u),\quad x\in  \mathbb{R}^{N},\\
&\|u\|_{2}=c,
\end{aligned}
\right.
\end{equation*}
where $a(x)\in C(\mathbb{R}^{N},\mathbb{R})$ satisfies $0<a_\infty=\lim \limits_{|x|\rightarrow \infty} a(x)\leq a(x)$ and $f\in C(\mathbb{R}, \mathbb{R})$ satisfies some Ambrosetti-Rabinowitz type condition.

Furthermore, by barycentric function and minimax method, Zhang and Zhang \cite{18} concerned with the existence and multiplicity of normalized bounded solutions for the following  Sch\"{o}dinger equations in exterior domain:
\begin{equation*}
\left\{
\begin{aligned}
&-\Delta u=\lambda u+|u|^{p-2}u\ &\mbox{in}\ \Omega ,\\
&u=0\ &\mbox{on}\ \partial\Omega,\\
&\int_{\Omega}|u|^{2}dx=c^{2},
\end{aligned}
\right.
\end{equation*}
where $p\in(2,2+4/N)$ and $\Omega\subset\mathbb{R}^{N}$ is an exterior domain.

Recently, the existence of normalized solutions for Schr\"{o}dinger equations with fractional Laplacian operator have attracted much attention, see \cite{11},\cite{14},\cite{16},\cite{b},\cite{c} and the references therein. In \cite{5}, Du and Tian et al. were devoted to study the existence and nonexistence of normalized solutions for the following equation
\begin{equation}\label{ASB}
\left\{
\begin{aligned}
&(-\Delta)^{s}u+V(x)u=\lambda u+\mu f(u)\ \ \text{in } \mathbb{R}^{N}, \\
&\int_{\mathbb{R}^{N}}|u|^{2}dx=c^{2},u\in H^{s}(\mathbb{R}^{N}),
\end{aligned}
\right.
\end{equation}
where $\mu>0$, $V$ is a trapping potential function and $f$ is a subcritical nonlinearity. If $f(u)=|u|^{p-2}u$ with $p\in(2+4s/N,2_{s}^{*})$ and $V(x)\in L^{\infty}(\mathbb{R}^{N})$ or $V(x)\in L^{\frac{N}{2}}(\mathbb{R}^{N})$ satisfies appropriate assumptions, Peng and Xia in \cite{16} proved existence of normalized bounded solutions for problem \eqref{ASB} with $\mu=1$.

In \cite{11}, the authors considered the existence of normalized ground state solution for the fractional Schr\"{o}dinger equations with combined nonlinearities:
$$
\left\{
\begin{aligned}
&(-\Delta)^{s}u=\lambda u+\mu|u|^{p-2}u+|u|^{q-2}u\ \ \text{in } \mathbb{R}^{N}, \\
&\int_{\mathbb{R}^{N}}|u|^{2}dx=c^{2},u\in H^{s}(\mathbb{R}^{N}).
\end{aligned}
\right.
$$
Li and Zou \cite{14}, Zhen and Zhang \cite{17} investigated existence of the normalized ground state solution for the critical fractional Schr\"{o}dinger equation, respectively.

Motivated by above papers, in this paper, we will investigate the existence of normalized bound solutions for problem \eqref{eq1} under $p\in(2,2+4s/N)$ and suitable assumption $(A_{1})$.

Normalized solutions of problem \eqref{eq1} are obtained by looking for critical point the following $C^1$ functional
\begin{equation}\label{eq5}
I(u)=\frac{1}{2}\int_{\mathbb{R}^{N}}|(-\Delta)^{s/2} u|^{2}dx-\frac{1}{p}\int_{\mathbb{R}^{N}}a(x)|u|^{p}dx
\end{equation}
constrained on the $L^2$-spheres in $H^s(\mathbb{R}^{N})$:
$$
S_{c}:=\{ u \in H^{s}(\mathbb{R}^{N}):\|u\|_{2}=c\},
$$
where $c>0$ is any fixed. Define
$$
m_{c}:=\inf _{u \in S_{c}}I(u),
$$
from $(A_{1})$ and \eqref{eq5}, we will prove $m_{c}>-\infty$ for any $p\in(2,2+4s/N)$ and $m_{c}=-\infty$ for any $p\in(2+4s/N, 2_s^*)$. $\bar{p}=2+4s/N$ is called the $L^2$-critical exponent.

From \cite{6}, let $p\in(2,2_{s}^{\ast})$, the following fractional Schr\"{o}dinger equation:
\begin{equation}\label{eq3.1}
(-\Delta)^{s}u+u=u^{p-1},\ x\in\mathbb{R}^{N},\quad u>0,
\end{equation}
has an unique solution $w\in H^{s}(\mathbb{R}^{N})$, which is radial, radially decreasing and positive. For $\lambda<0$, let
$w_{\lambda}(x):=(-\lambda)^{\frac{1}{p-2}} w((-\lambda)^{\frac{1}{2s}}x)$,  $w_{\lambda}$ is the unique positive solution (up to a translation) of
\begin{equation}\label{eq3}
(-\Delta)^{s}u=\lambda u+u^{p-1}  \text { in } \mathbb{R}^{N}, \quad u>0.
\end{equation}

To restore compactness of $(PS)$-sequence for the functional $I$, we need to consider the uniqueness of solutions for the limited problem corresponding to problem \eqref{eq1}:
\begin{equation} \label{eq4}
\left\{
\begin{aligned}
&(-\Delta)^{s}u=\lambda u+|u|^{p-2}u,\\
&\int_{\mathbb{R}^{N}}|u|^{2}dx=c^{2},\ u\in H^{s}(\mathbb{R}^{N}).
\end{aligned}
\right.
\end{equation}
Define the functional $I_{\infty}: H^{s}(\mathbb{R}^{N})\to \mathbb{R}$ as follows:
\begin{equation}\label{eq55}
I_{\infty}(u)=\frac{1}{2}\int_{\mathbb{R}^{N}}|(-\Delta)^{s/2} u|^{2}dx-\frac{1}{p}\int_{\mathbb{R}^{N}}|u|^{p}dx
\end{equation}
and
\begin{equation}\label{eq7}
m^{\infty}_{c}:=\inf _{u \in S_{c}}I_{\infty}(u).
\end{equation}

Now we state our main results of this paper.\\

\textbf{Theorem 1.1.}\ Assume that $(A_{1})$ holds, we have $m_c=m^{\infty}_{c}$, and $m_{c}$ is not attained for the functional $I$ in $S_{c}$.\\

\textbf{Theorem 1.2.}\ Assume that $(A_{1})$ and

$(A_{2})$\ $\sup\limits_{x\in\mathbb{R}^{N}}|a_{\infty}-a(x)|$ is sufficiently small.

\noindent Then, problem \eqref{eq1} has a positive solution $u \in S_{c}$ for some $\lambda<0$. \\

\textbf{Remark.}\ In\cite{LMS}, Lehrer, Maia and Squassina exploited Poho\v{z}aev manifold and a variational technique to prove existence of bound state solutions for the asymptotically linear fractional Schrodinger equations. Recently, the existence of normalized solutions of mass subcritical fractional Schr\"{o}dinger equations were obtained by Yu, Tang and Zhang (see\cite{b}) in exterior domains. In this paper, we consider the normalized bound state solutions of fractional Schr\"{o}dinger equations in the whole space. Hence our situation is different from \cite{LMS} and \cite{b}, and the main idea of our argument comes from \cite{18}, and to obtain the non-negative Palais-Smale sequence of the functional $I$ in $S_c$
, we will apply the argument similar with the one of \cite{AS}.

\section{the Proof of Theorem 1.1}


\ \ \ \ \ \ In this section, we first consider the minimization problem \eqref{eq7}, and then prove Theorem 1.1. First of all, by the change of variable, we can assume $a_\infty=1$ for convenience.\\

\textbf{Lemma 2.1.}\ 
For any $c>0$, we have $m^{\infty}_{c}\in(-\infty,0)$. Moreover, if $\{u_{n}\} \subset S_{c}$ is a non-negative minimizing sequence for the functional $I_{\infty}$ in $S_{c}$, then up to a subsequence, there exists $\{y_{n}\} \subset \mathbb{R}^{N}$ such that
$$
u_{n}(\cdot-y_{n}) \rightarrow w_{\lambda}\ \text { strongly in } H^{s}(\mathbb{R}^{N})
$$
as $n \rightarrow \infty$, where $\lambda<0$ is determined by
\begin{equation}\label{eq22}
c^{2}=\|w_\lambda\|_{2}^{2}=(-\lambda)^{s_{p}}\|w\|_{2}^{2}, \quad s_{p}:=\frac{2}{p-2}-\frac{N}{2s}.
\end{equation}
In particular, $m^{\infty}_{c}$ is attained by the function $w_{\lambda}$ and can be expressed as
\begin{equation}\label{eq23}
m^{\infty}_{c}=\left(\frac{c^{2}}{\|w\|_{2}^{2}}\right)^{\frac{s_{p}+1}{s_{p}}} I_{\infty}(w).
\end{equation}
where $w$ is an unique solution for problem \eqref{eq3.1}.

\textbf{Proof.}\ By \eqref{eq6}, for any $u \in S_{c}$, we have
\begin{eqnarray}\label{eq24}
 I_{\infty}(u) &=& \frac{1}{2}\int_{\mathbb{R}^{N}}|(-\Delta)^{s/2} u|^{2}dx-\frac{1}{p}\|u\|_{p}^{p}\nonumber \\
   &\geq& \frac{1}{2}\int_{\mathbb{R}^{N}}|(-\Delta)^{s/2} u|^{2}dx-\frac{C_{N,p,s}}{p}\left(\int_{\mathbb{R}^{N}}|(-\Delta)^{s/2} u|^{2}dx\right)^{\frac{N(p-2)}{4s}}\|u\|_{2}^{p-\frac{N(p-2)}{2s}}.
\end{eqnarray}
Since $2<p<2+\frac{4s}{N}$, then $0<\frac{N(p-2)}{4s}<1$, and it follows from \eqref{eq24} that $I_{\infty}$ is coercive in $S_{c}$. Therefore, $m^{\infty}_{c}>-\infty$. Moreover, for any $u \in S_{c}$, $t>0$, one has $t^{\frac{N}{2}}u(tx)\in S_{c}$ and
$$
I_{\infty}(t^{\frac{N}{2}}u(tx))=\frac{1}{2}t^{2s}\int_{\mathbb{R}^{N}}|(-\Delta)^{s/2} u|^{2}dx-\frac{1}{p}t^{N(\frac{p}{2}-1)}\|u\|_{p}^{p} \rightarrow 0^{-}\ \mbox{as}\ t \to 0^+.
$$
One can get that $m^{\infty}_{c}<0$, that is, $m^{\infty}_{c} \in(-\infty, 0)$.

Let $c>0$, $\theta>1$ and $\{v_{n}\} \subset S_{c}$ be a minimizing sequence for $m^{\infty}_{c}$. From $\theta>1$, $p>2$, we have
$$
m^{\infty}_{\theta c}\leq I_{\infty }(\theta v_{n})=\frac{\theta^{2}}{2}\int_{\mathbb{R}^{N}}|(-\Delta)^{s/2} v_n|^{2}dx-\frac{\theta^{p}}{p}\|v_{n}\|_{p}^{p}<\theta^{2} I_{\infty}(v_{n}),
$$
hence, it follows that $m^{\infty}_{\theta c}\leq \theta^{2}m^{\infty}_{c}$, and the equality holds if and only if $\lim _{n \rightarrow \infty}\|v_{n}\|_{p}^{p}=0$. We claim that $m^{\infty}_{\theta c}<\theta^{2}m^{\infty}_{c}$. Suppose, by contradiction, that $m^{\infty}_{\theta c}=\theta^{2}m^{\infty}_{c}$, then $\lim _{n \rightarrow \infty}\|v_{n}\|_{p}^{p}=0$. Hence, we have
$$
m^{\infty}_{c}=\lim _{n \rightarrow \infty} I_{\infty }(v_{n}) \geq \liminf _{n \rightarrow \infty} \frac{1}{2}\int_{\mathbb{R}^{N}}|(-\Delta)^{s/2}v_n|^{2}dx \geq 0,
$$
which contradicts with $m^{\infty}_{c}<0$.

Let $c_{1}$, $c_{2}>0$ be such that $c_{1}^{2}+c_{2}^{2}=c^{2}$ and let  $\theta_{j}=\frac{c}{c_{j}}>1$ for $j=1,\ 2$ we have $m^{\infty}_{c}=m^{\infty}_{\theta_{j} c_{j}}<\theta_{j}^{2} m^{\infty}_{c_{j}}$. From  $c_{1}^{2}+c_{2}^{2}=c^{2}$, we conclude
$$
m^{\infty}_{c}=\left(\frac{1}{\theta_{1}^{2}}+\frac{1}{\theta_{2}^{2}}\right)m^{\infty}_{c}<m^{\infty}_{c_{1}}+m^{\infty}_{c_{2}}.
$$

By Theorem $1.2$ of \cite{13}, up to a subsequence, there exist $u \in S_{c}$ and $\{y_{n}^{\prime}\} \subset \mathbb{R}^{N}$ such that
$u_{n}(\cdot-y_{n}^{\prime}) \rightarrow u$ strongly in $H^{s}(\mathbb{R}^{N})$ as $n \rightarrow \infty$.
In particular, $I_{\infty}(u)=m^{\infty}_{c}$. Since $u_{n} \geq 0$ in $\mathbb{R}^{N}$ for any $n \geq 1$, one has $u \geq 0$ in $\mathbb{R}^{N}$. By the Lagrange multiplier rule, there exists $\lambda \in \mathbb{R}$ such that $u$ is a non-negative nontrivial weak solution of
$$
(-\Delta)^{s}u-|u|^{p-2}u=\lambda u \ \text { in }\ H^{s}(\mathbb{R}^{N}),
$$
which implies that
$$
\lambda c^{2}=\int_{\mathbb{R}^{N}}|(-\Delta)^{s/2} u|^{2}dx-\|u\|_{p}^{p}=pm^{\infty}_{c}+\frac{2-p}{2}\int_{\mathbb{R}^{N}}|(-\Delta)^{s/2} u|^{2}dx<0,
$$
hence, we have $\lambda<0$. By the strong maximum principle, $u>0$ in $\mathbb{R}^{N}$. The uniqueness of solutions for problem \eqref{eq3} guarantees that there exists $y \in \mathbb{R}^{N}$ such that $u(\cdot-y)=w_{\lambda}$. Denote $y_{n}=y_{n}^{\prime}+y$ for any $n \geq 1$, we have $u_{n}(\cdot-y_{n}) \rightarrow w_{\lambda}$ strongly in $H^{s}(\mathbb{R}^{N})$ as $n \rightarrow \infty$.
Particularly, $w_{\lambda}\in S_{c}$ and $I_{\infty}(w_{\lambda})=m^{\infty}_{c}$. By the definition of $w_{\lambda}$, \eqref{eq22} and \eqref{eq23} hold.\ \ \ $\Box$\\

Define the map $\Psi: \mathbb{R}^{N} \rightarrow S_{c}$ by
$$
\Psi(y)=w_{\lambda}(\cdot-y)\ \text{for all}\ y \in \mathbb{R}^{N},
$$
where $w_{\lambda}(x)$ is given in Lemma 2.1. \\ 

\textbf{Proof of Theorem 1.1.}\ Since $a(x)\leq 1$, one has $m_{c} \geq m^{\infty}_{c}$. On the other hand, from
$$
\int_{\mathbb{R}^{N}}a(x)|w_{\lambda}(x-y)|^pdx\to\int_{\mathbb{R}^{N}}|w_{\lambda}(x)|^pdx\ \ \text{as}\ \ |y|\to\infty,
$$
it follows that
$$
m_{c}\leq I(\Psi(y))\rightarrow I_{\infty}(w_{\lambda})=m^{\infty}_{c}\ \ \text{as}\ \ |y|\to\infty.
$$
Hence, we have $m_{c}=m^{\infty}_{c}\in(-\infty, 0)$.

Suppose, by contradiction, that there is $u\in S_{c}$ such that $I(u)=m_{c}$. Note that $|u|\in S_{c}$ and
\begin{eqnarray}\label{AS}
   \int_{\mathbb{R}^{N}}|(-\Delta)^{\frac{s}{2}} |u||^{2}dx&=& \int_{\mathbb{R}^{N}}\int_{\mathbb{R}^{N}} \frac{||u(x)|-|u(y)||^{2}}{|x-y|^{N+2 s}}dxdy \nonumber \\
  &\leq& \int_{\mathbb{R}^{N}}\int_{\mathbb{R}^{N}} \frac{|u(x)-u(y)|^{2}}{|x-y|^{N+2 s}}dxdy =\int_{\mathbb{R}^{N}}|(-\Delta)^{\frac{s}{2}} u|^{2}dx,
\end{eqnarray}
from $a(x)\leq 1$, we have
$$
m_{c}=I(u)\geq I_{\infty}(u)\geq I_{\infty}(|u|)\ge m_{c}^{\infty}.
$$
Hence from Lemma 2.1, we have $u=w_\lambda$ up to a translation, where $\lambda<0$ is determined by \eqref{eq22}.  From $I(w_\lambda)=I_{\infty}(w_\lambda)=m_{c}^{\infty}$, we have
$$
\int_{\mathbb{R}^{N}}(1-a(x))|w_\lambda|^pdx=0,
$$
which is a contradiction with $w_\lambda>0$, $a(x)\leq 1$ and $a(x)\not\equiv 1$ in $\mathbb{R}^{N}$. Hence $m_{c}$ is not attained for the functional $I$ in $S_{c}$.\ \ \ $\Box$

\section{the Proof of Theorem 1.2}

\ \ \ \ \ \ In the subsection, we will introduce the Splitting Lemma, which can be found in \cite{LMS} and is crucial to ensure the convergence of the some non-negative Palais-Smale sequence of the functional $I$ in $S_c$. And then, we will prove some useful lemmas and Theorem 1.2.

First of all, we define the energy functionals associated with problem \eqref{eq1} and problem \eqref{eq3} respectively:
$$
I_{\lambda}(u)=\frac{1}{2}\int_{\mathbb{R}^{N}}|(-\Delta)^{s/2} u|^{2}dx-\frac{1}{p}\int_{\mathbb{R}^{N}}a(x)|u|^{p}dx-\frac{\lambda}{2}\int_{\mathbb{R}^{N}}|u|^{2}dx,
$$
$$
I_{\lambda,\infty}(u)=\frac{1}{2}\int_{\mathbb{R}^{N}}|(-\Delta)^{s/2} u|^{2}dx-\frac{1}{p}\int_{\mathbb{R}^{N}}|u|^{p}dx-\frac{\lambda}{2}\int_{\mathbb{R}^{N}}|u|^{2}dx
$$
for any $u\in H^{s}(\mathbb{R}^{N})$.\\

\textbf{Lemma 3.1.}\ For any $\lambda<0$, let $\{u_{n}\}\subset H^{s}(\mathbb{R}^{N})$ be a non-negative Palais-Smale sequence of the functional $I_{\lambda}$. Then up to a subsequence, there exist an integer $k \geq 0$, a non-negative function $u_{0} \in H^{s}(\mathbb{R}^{N})$, $k$ sequences $\{x_{n}^{i}\} \subset \mathbb{R}^{N}$ for $1 \leq i \leq k$ such that $|x_{n}^{i}| \rightarrow\infty$ as $n \rightarrow \infty$ and
$$
u_{n}=u_{0}+\sum_{i=1}^{k} w_{\lambda}(\cdot+x_{n}^{i})+o(1)\ \text { strongly in } H^{s}(\mathbb{R}^{N}).
$$
Furthermore, $u_{0}$ is a weak solution of
$$
(-\Delta)^{s}u=\lambda u+a(x)|u|^{p-2}u\ \text{ in }\mathbb{R}^{N},
$$
and $w_\lambda(x)$ is a weak solution of
$$
(-\Delta)^{s}u=\lambda u+|u|^{p-2}u\ \text{ in }\mathbb{R}^{N}.
$$
Moreover, we have
\begin{equation}\label{eq21}
\|u_{n}\|_{2}^{2}=\|u_{0}\|_{2}^{2}+k\|w_{\lambda}\|_{2}^{2}+o(1),
\end{equation}
\begin{equation}\label{eq21.1}
I_{\lambda}(u_{n})=I_{\lambda}(u_{0})+k I_{\lambda,\infty}(w_{\lambda})+o(1).
\end{equation}
\vspace{0.1cm}

\textbf{Lemma 3.2.}\  Let $N>2s$, $2<p<2+\frac{4s}{N}$, $c>0$, there exists a positive constant $\eta=\eta(c)\in(2^{-1/s_{p}}, 1)$ such that if $\{u_{n}\} \subset S_{c}$ is a non-negative Palais-Smale sequence of $I|_{S_{c}}$ at level $d$ with $m^{\infty}_{c}<d<\eta m^{\infty}_{c}$, then up to a subsequence, there exists $u_{0} \in S_{c}$ such that
$u_{n} \rightarrow u_{0}$ strongly in $H^{s}(\mathbb{R}^{N})$ as $n \rightarrow \infty$. Moreover, $u_{0}$ is a positive solution of problem \eqref{eq1} for some $\lambda<0$.

\textbf{Proof.}\ Since $2<p<2+\frac{4s}{N}$, we have $s_{p}>0$ and $0<2^{-1 / s_{p}}<1$. Due to $m^{\infty}_{c}<0$, one has
$$
m^{\infty}_{c}<2^{-1/s_{p}}m^{\infty}_{c}<0.
$$
Set $\{u_{n}\} \subset S_{c}$ be a non-negative Palais-Smale sequence of $I|_{S_{c}}$ at level $d \in(m^{\infty}_{c}, 2^{-1/s_{p}}m^{\infty}_{c})$. Because of $I(u_{n})\rightarrow d<0$ as $n \rightarrow \infty$ and \eqref{eq24}, $\{u_{n}\}$ is bounded in $H^{s}(\mathbb{R}^{N})$. Note that $(I|_{S_c})^{\prime}(u_{n})=o(1)$, according to the Lagrange multiplier rule, there exists $\{\lambda_{n}\} \subset \mathbb{R}$ such that
\begin{equation}\label{eq25}
\int_{\mathbb{R}^{N}}(-\Delta)^{\frac{s}{2}}u_{n}(-\Delta)^{\frac{s}{2}}\psi d x-\lambda_{n} \int_{\mathbb{R}^{N}} u_{n} \psi d x-\int_{\mathbb{R}^{N}}a(x)u_{n}^{p-1} \psi d x=o(1)\|\psi\|
\end{equation}
for any $\psi \in H^{s}(\mathbb{R}^{N})$. Since $\{u_{n}\} \subset S_{c}$ is bounded in $H^{s}(\mathbb{R}^{N})$, using \eqref{eq25}, we have
$$
-\lambda_{n}c^{2}=\int_{\mathbb{R}^{N}}a(x)|u_{n}|^{p}dx-\int_{\mathbb{R}^{N}}|(-\Delta)^{s/2} u_n|^{2}dx +o(1),
$$
which implies that $\{\lambda_{n}\}$ is bounded in $\mathbb{R}$. Thus, up to a subsequence, there exists $\lambda \in \mathbb{R}$ such that $\lambda_{n} \rightarrow \lambda$ as $n \rightarrow \infty$. Moreover, since $p>2$ and $I(u_{n}) \rightarrow d<0$ as $n \rightarrow \infty$, we obtain
$$
\begin{aligned}
-\lambda_{n} c^{2}&=\int_{\mathbb{R}^{N}}a(x)|u_{n}|^{p}dx-\int_{\mathbb{R}^{N}}|(-\Delta)^{s/2} u_n|^{2}dx+o(1)\\
&=-pd+\frac{p-2}{2}\int_{\mathbb{R}^{N}}|(-\Delta)^{s/2} u_n|^{2}dx+o(1)\\
&\geq-pd+o(1).\\
\end{aligned}
$$
Letting $n \rightarrow \infty$, we can get
\begin{equation}\label{eq26}
-\lambda \geq \frac{-pd}{c^{2}}>0,
\end{equation}
That is, $\lambda<0$. From \eqref{eq25} and the boundedness of $\{u_{n}\}$ in  $H^{s}(\mathbb{R}^{N})$, we deduce that $\{u_{n}\}$ is a non-negative Palais-Smale sequence of the functional $I_{\lambda}$. By Lemma $3.1$, up to a subsequence, there exists an integer $k \geq 0$, a non-negative function $u_{0} \in H^{s}(\mathbb{R}^{N})$ and $k$ sequences  $\{x_{n}^{i}\} \subset \mathbb{R}^{n}$ for $1 \leq i \leq k$ such that  $|x_{n}^{i}| \rightarrow\infty$ as $n \rightarrow \infty$ and
\begin{equation}\label{eq27}
u_{n}=u_{0}+\sum_{i=1}^{k} w_{\lambda}(\cdot+x_{n}^{i})+o(1)\ \text{ strongly in } H^{s}(\mathbb{R}^{N}).
\end{equation}
Moreover, $u_{0}$ is a weak solution of
\begin{equation}\label{eq28}
(-\Delta)^s u=\lambda u+a(x)u^{p-2}u \ \text { in } \mathbb{R}^{N}
\end{equation}
and
\begin{equation}\label{eq29}
c^{2}=\|u_{0}\|_{2}^{2}+k\|w_{\lambda}\|_{2}^{2},
\end{equation}
\begin{equation}\label{eq210}
I_{\lambda}(u_{n})=I_{\lambda}(u_{0})+k I_{\lambda,\infty}(w_{\lambda})+o(1).
\end{equation}

Because
$$I_{\lambda}(v)=I(v)-\frac{\lambda}{2}\|v\|_{2}^{2},\ \ \ I_{\lambda,\infty}(v)=I_{\infty}(v)- \frac{\lambda}{2}\|v\|_{2}^{2}
$$
for any $v \in H^{s}(\mathbb{R}^{N})$, by \eqref{eq29}, \eqref{eq210} and Lemma $2.1$, we can see
\begin{eqnarray}\label{eq211}
 d&=& I(u_{0})+km_{\|w_{\lambda}\|_{2}}^{\infty}\geq m_{\|u_{0}\|_{2}}^{\infty }+km_{\|w_{\lambda}\|_{2}}^{\infty}\nonumber \\
 \ &=&C_{0}\left((c^{2}-k\|w_{\lambda}\|_{2}^{2})^{\frac{s_{p}+1}{s_{p}}}
+k(\|w_{\lambda}\|_{2}^{2})^{\frac{s_{p}+1}{s_{p}}}\right),
\end{eqnarray}
where
$$
C_{0}:=I_{\infty}(w)(\|w\|_{2}^{2})^{-\frac{s_{p}+1}{s_{p}}}<0.
$$

We claim that there exists $\eta=\eta(c) \in(2^{-1 / s_{p}}, 1)$ such that if  $m^{\infty}_{c}<d<\eta m^{\infty}_{c}$, then $k=0$. Once this claim is proved, using \eqref{eq27}, we can see that $u_{n} \rightarrow u_{0}$ strongly in $H^{s}(\mathbb{R}^{N})$. It implies that $u_{0} \in S_{c}$ and $u_{0}$ is a nontrivial non-negative solution of problem \eqref{eq1} with $\lambda<0$. The strong maximum principle then yields that $u_{0}>0$ in $\mathbb{R}^{N}$. In the following, we prove the above claim by three steps.

\textbf{Step 1}. If $m^{\infty}_{c}<d<2^{-1 / s_{p}}m^{\infty}_{c}$, then $\|w_{\lambda}\|_{2}^{2} \geq l c^{2}$, where
$$
l:=\frac{1}{2}\left[\frac{-p I_{\infty }(w)}{\|w\|_{2}^{2}}\right]^{s_{p}} \in(0, \frac{1}{2}).
$$
From Lemma 2.1 and \eqref{eq26}, we obtain
$$
\begin{aligned}
\|w_{\lambda}\|_{2}^{2}
&=(-\lambda)^{s_{p}}\|w\|_{2}^{2} \geq(\frac{-pd}{c^{2}})^{s_{p}}\|w\|_{2}^{2} \\
&\geq\left(\frac{-p 2^{-1 /s_{p}}m^{\infty}_{c}}{c^{2}}\right)^{s_{p}}\|w\|_{2}^{2} \\
& =lc^{2}.
\end{aligned}
$$
Observe that $w$ is a solution of \eqref{eq3.1}, it is clear that
$$
\int_{\mathbb{R}^{N}}|(-\Delta)^{s/2} w|^{2}dx+\|w\|_{2}^{2}=\|w\|_{p}^{p}.
$$
Hence
$$
I_{-1,\infty}(w)=\left(\frac{1}{2}-\frac{1}{p}\right)\left(\int_{\mathbb{R}^{N}}|(-\Delta)^{s/2} w|^{2}dx+\|w\|_{2}^{2}\right)>\left(\frac{1}{2}-\frac{1}{p}\right)\|w\|_{2}^{2}.
$$
In view of $I_{\infty}(w)=I_{-1,\infty}(w)-\frac{1}{2}\|w\|_{2}^{2}<0$, we conclude that
$$
0<\frac{-p I_{\infty }(w)}{\|w\|_{2}^{2}}<1.
$$
Since $s_{p}>0$, we get $l \in\left(0, \frac{1}{2}\right)$.

\textbf{Step 2}. Define
$$
l_{1}:=l^{(s_{p}+1)/s_{p}}+(1-l)^{(s_{p}+1)/s_{p}}\in(2^{-1/s_{p}}, 1) .
$$
If $m^{\infty}_{c}<d<l_{1}m^{\infty}_{c}$, then $k \leq 1$.
Actually, for $k \geq 1$, we set
$$
f_{k}(t):=(c^{2}-k t)^{\frac{s_{p}+1}{s_{p}}}+k t^{\frac{s_{p}+1}{s_{p}}}\ \mbox{for all}\ 0 \leq t \leq \frac{c^{2}}{k} .
$$
Then $f_{k} \in C^{1}\left(\left[0, \frac{c^{2}}{k}\right]\right)$ and
$$
f_{k}^{\prime}(t)=\frac{s_{p}+1}{s_{p}} k[t^{1 / s_{p}}-(c^{2}-kt)^{1 / s_{p}}], \quad \forall t \in\left[0, \frac{c^{2}}{k}\right] .
$$
Hence, $f_{k}(t)$ has a critical point $t_{k}:=\frac{c^{2}}{k+1}$, and is strictly decreasing in $[0, t_{k}]$ and strictly increasing in $\left[t_{k}, \frac{c^{2}}{k}\right]$. Moreover, we have
$$
\min _{t \in[0, c^{2}/k]} f_{k}(t)=f_{k}(t_{k})=(k+1)^{-\frac{1}{s_{p}}}(c^{2})^{\frac{s_{p}+1}{s_{p}}}>0 .$$
Step $1$ and \eqref{eq29} leads to $k \leq l^{-1}$. Define
$$
h(t):=(1-lt)^{(s_{p}+1)/s_{p}}+t l^{(s_{p}+1)/ s_{p}}, \quad \forall t \in[1, l^{-1}],
$$
it follows that $h \in C^{1}([1, l^{-1}])$ and
$$
h^{\prime}(t)=l^{(s_{p}+1) / s_{p}}-l \frac{s_{p}+1}{s_{p}}(1-lt)^{\frac{1}{s_{p}}}, \quad \forall t\in[1, l^{-1}].
$$
Hence, $h(t)$ has a critical point $t_{0}:=\frac{1}{l}-\left(\frac{s_{p}}{s_{p}+1}\right)^{s_{\mathrm{p}}} \in(1, l^{-1})$, and is strictly decreasing in $[1, t_{0}]$ and strictly increasing in  $[t_{0}, l^{-1}]$.

If $k \geq 2$, by \eqref{eq211}, we deduce that
$$
\begin{aligned}
d & \geq C_{0} \max\{f_{k}(lc^{2}), f_{k}(c^{2}/k)\}=\max\{h(k), k^{-1 / s_{p}}\}m^{\infty}_{c}\\
& \geq \max\{h(1), h(l^{-1}), 2^{-1 / s_{p}}\}m^{\infty}_{c} \\
& =l_{1}m^{\infty}_{c},
\end{aligned}
$$
which is a contradiction to $d<l_{1}m^{\infty}_{c}$. Therefore, we have $k \leq 1$.

\textbf{Step 3}. There exists $\eta=\eta(c) \in[l_{1}, 1)$ such that if $m^{\infty}_{c}<d<\eta m^{\infty}_{c}$, then $k=0$.

From Step $2$, we have $k \leq 1$. If $u_{0} \equiv 0$, then by \eqref{eq29} and \eqref{eq211}, we get $d=0$ or $d=m^{\infty}_{c}$, which contradicts to $0>d>m^{\infty}_{c}$. Hence, $u_{0} \not \equiv 0$. By \eqref{eq28}, we get
\begin{equation}\label{eq212}
\int_{\mathbb{R}^{N}}|(-\Delta)^{s/2} u_0|^{2}dx+(-\lambda)\|u_{0}\|_{2}^{2}=\int_{\mathbb{R}^{N}}a(x)|u_{0}|^{p}dx.
\end{equation}
On the one hand, from $a(x)\leq 1$, \eqref{eqA2} and \eqref{eq212}, there exists $C_{1}>0$ independent of $u_{0}$ and $\lambda$ such that
$$
\min \{1,-\lambda\}\|u_{0}\|^2 \leq \|u_{0}\|_{p}^{p}\leq C_{1}\|u_{0}\|^p.
$$
Because $u_{0} \not \equiv 0$, $\lambda<0$ and $p>2$, we conclude that
\begin{equation}\label{eq213}
\int_{\mathbb{R}^{N}}|(-\Delta)^{s/2} u_0|^{2}dx+\|u_{0}\|_{2}^{2}=\|u_{0}\|^2\geq\left(\frac{\min \{1,-\lambda\}}{C_{1}}\right)^{\frac{2}{p-2}}>0.
\end{equation}

On the other hand, by \eqref{eq6}, \eqref{eq212} and $\lambda<0$, we can see
$$
\int_{\mathbb{R}^{N}}|(-\Delta)^{s/2} u_0|^{2}dx\leq\|u_{0}\|_{p}^{p} \leq C_{N,p,s}\left(\int_{\mathbb{R}^{N}}|(-\Delta)^{s/2} u_0|^{2}dx\right)^{\frac{N(p-2)}{4s}}\|u_0\|_{2}^{p-\frac{N(p-2)}{2s}}.
$$
From $2<p<2+\frac{4s}{N}$, let $\gamma_{p}:=\frac{N(p-2)}{4s}$, we get
\begin{equation}\label{eq214}
\int_{\mathbb{R}^{N}}|(-\Delta)^{s/2} u_0|^{2}dx\leq C_{N, p,s}^{\frac{1}{1-\gamma_{p}}}\|u_{0}\|_{2}^{\frac{p-2\gamma_{p}}{1-\gamma_{p}}}.
\end{equation}
Owing to \eqref{eq26}, \eqref{eq213}, \eqref{eq214} and $d<2^{-1 / s_{p}}m^{\infty}_{c}$  that
$$
\begin{array}{l}
C_{N,p,s}^{\frac{1}{1-\gamma_{p}}}(\|u_{0}\|_{2}^2)^{\frac{p-2\gamma_{p}}{2-2\gamma_{p}}}
+\|u_{0}\|_{2}^{2}\geq\left(\frac{1}{C_{1}} \min \left\{1,-p 2^{-1 / s_{p}} \frac{m^{\infty}_{c}}{c^{2}}\right\}\right)^{\frac{2}{p-2}}>0,
\end{array}
$$
which means that there exists $C(c)>0$ such that
$$
\|u_{0}\|_{2}^{2} \geq C(c).
$$

Set
$$
\eta(c):=\left\{\begin{array}{ll}
l_{1}, & \text { if } c^{2} \leq C(c), \\
\max \left\{l_{1}, \frac{f_{1}(c^{2}-C(c))}{(c^{2})^{\frac{s_{p}+1}{s_{p}}}}\right\}, & \text { if } c^{2}>C(c).
\end{array}\right.
$$
If $c^{2} \leq C(c)$, then $\|u_{0}\|_{2}^{2} \geq C(c) \geq c^{2}$. By \eqref{eq29} we get $\|u_{0}\|_{2}^{2}=c^{2}$ and hence  $k=0$.
If $c^{2}>C(c)$, then $0<c^{2}-C(c)<c^{2}$. Hence, we have
$$
0<f_{1}(c^{2}-C(c))<\max \{f_{1}(0), f_{1}(c^{2})\}=(c^{2})^{\frac{s_{p}+1}{s_{p}}},
$$
which implies that  $0<\frac{f_{1}(c^{2}-C(c))}{(c^{2})^{\frac{s_{p}+1}{s_{p}}}}<1$. In particular, $\eta(c) \in[l_{1}, 1)$. Assume that $k=1$, it follows from Step $1$ and \eqref{eq29} that
$$
l c^{2} \leq\|w_{\lambda}\|_{2}^{2} \leq c^{2}-C(c).
$$
By \eqref{eq211}, one has
$$
d\geq C_{0} \max \{f_{1}(l c^{2}), f_{1}(c^{2}-C(c))\} \geq \eta(c)m^{\infty}_{c},
$$
which contradicts to $d<\eta(c)m^{\infty}_{c}$. Hence $k=0$.\ \ \ $\Box$\\

Let the function $w_\lambda(x)$ and $\eta=\eta(c)\in(0,1)$ be given in Lemma 2.1 and Lemma 3.1, respectively, now we replace $(A_{2})$ by a accurate condition, and we have the following lemma.\\

\textbf{Lemma 3.3.}\ Assume that
$$
\sup_{x\in\mathbb{R}^{N}}|1-a(x)|<\frac{(\eta-1)m^{\infty}_{c}}{\| w_{\lambda}\|_{p}^{p}},
$$
we have $I(\Psi(y))<\eta m^{\infty}_{c}$.

\textbf{Proof.}\ Since $I_{\infty}$ is translation invariant, we have $I_{\infty}(\Psi(y))=I_{\infty}(w_{\lambda})=m^{\infty}_{c}$, and
$$
\begin{aligned}
I(\Psi(y))&=I_{\infty}(\Psi(y))+I(\Psi(y))-I_{\infty}(\Psi(y))\\
&\le m^{\infty}_{c}+\int_{\mathbb{R}^{N}}(1- a(x))|\Psi(y)|^{p}dx\\
&\le m^{\infty}_{c}+\sup_{\mathbb{R}^{N}}|1-a(x)|\|\Psi(y)\|_{p}^{p}\\
&<m^{\infty}_{c}+\frac{(\eta-1)m^{\infty}_{c}}{\| w_{\lambda}\|_{p}^{p}}\|w_{\lambda}\|_{p}^{p}\\
&=\eta m^{\infty}_{c}.\ \ \ \ \ \ \ \ \Box
\end{aligned}
$$

\textbf{Definition 3.4.}(\cite{15})\ Define the barycenter function of a given function $u\in H^{s}(\mathbb{R}^{N})\setminus\{0\}$ as follows: let
$$
\mu(u)(x)=\frac{1}{|B_{1}|} \int_{B_{1}(x)}|u(y)|dy,
$$
we have $\mu(u) \in L^{\infty}(\mathbb{R}^{N})$ and $\mu$ is a continuous function. Subsequently, let
$$
\hat{u}(x)=\left[\mu(u)(x)-\frac{1}{2} \max \mu(u)\right]^{+},
$$
it follows that $\hat{u} \in C_{0}(\mathbb{R}^{N})$. Now define the barycenter of $u$ by
$$
\beta(u)=\frac{1}{\|\hat{u}\|_{1}} \int_{\mathbb{R}^{N}} x \hat{u}(x) d x \in \mathbb{R}^{N}.
$$
Since $\hat{u}$ has compact support, by the definition, $\beta(u)$ is well defined.
The function $\beta$ satisfies the following properties:

$(a)$ $\beta$ is a continuous function in $H^{s}(\mathbb{R}^{N}) \backslash\{0\}$.

$(b)$ If $u$ is radial, then $\beta(u)=0$.

$(c)$ For any $y \in \mathbb{R}^{N}$ and $u\in H^{s}(\mathbb{R}^{N})\setminus\{0\}$, we have  $\beta(u(x-y))=\beta(u(x))+y$.\\

Let
$$
P=\{u\in H^s(\mathbb{R}^{N}):u\geq 0\},
$$
now we define
$$
b_{c}:=\inf \{I(u): u\in S_{c}\cap P,\ \beta(u)=0\}.
$$
It is clear that $b_{c} \geqslant m_{c}^{\infty}$. Moreover, we have the following lemma.\\

\textbf{Lemma 3.5.}\  $b_{c}>m_{c}^{\infty}$.

\textbf{Proof.}\ By the definition of $b_{c}$, it is clear that $b_{c} \geqslant m_{c}^{\infty}$. Suppose, by contradiction, that $b_{c}=m_{c}^{\infty}$. There exists a sequence $\{u_{n}\}\subset S_{c}\cap P$ such that
$\beta(u_{n})=0$ for any $n \geq 1$ and $I(u_{n})\rightarrow m_{c}^{\infty}$ as $n \rightarrow \infty$. From $a(x)\leq 1$, we have $I_\infty(u_{n})\rightarrow m_{c}^{\infty}$ as $n \rightarrow \infty$. By Lemma 2.1, up to a subsequence, there exists a sequence $\{x_{n}\}\subset\mathbb{R}^{N}$ such that
\begin{equation}\label{eq215}
u_{n}(\cdot-x_{n})\to w_{\lambda}\ \mbox{ strongly in }H^s(\mathbb{R}^{N}),
\end{equation}
where $\lambda<0$ is determined by \eqref{eq22}.

Since $w_{\lambda}$ is radial, by the properties of $\beta$, we have
$$
x_{n}=\beta(u_{n})+x_{n}=\beta(u_{n}(\cdot-x_{n}))\to \beta(w_{\lambda})=0 \ \mbox{as}\ n\to\infty.
$$
Combining with \eqref{eq215}, it follows that
$$
u_{n}\to w_{\lambda}\ \mbox{ strongly in }H^s(\mathbb{R}^{N}).
$$
Since $I(u_{n})\rightarrow m_{c}=m_{c}^{\infty}$ as $n \rightarrow \infty$, we have
$$
I(w_{\lambda})= m_{c},
$$
which is a contradiction with that $m_{c}$ can not be attained. Therefore, $b_{c}>m_{c}^{\infty}$.\ \ $\Box$\\

\textbf{Proof of Theorem 1.2.}\  Under assumption $(A_{1})$, for all $u \in H^{s}(\mathbb{R}^{N}) \backslash\{0\}$, we get $I_{\infty}(u)<I(u)$. In particular,  $I_{\infty}(\Psi(y))<I(\Psi(y))$ for any $y \in \mathbb{R}^{N}$. It follows from Lemma 3.5 and $I(\Psi(y)) \rightarrow m^{\infty}_{c}$ as $|y| \rightarrow \infty$, there exists $\bar{\rho}>0$ such that
\begin{equation}\label{IP}
m^{\infty}_{c}<\max\limits_{|y|=\bar{\rho}} I(\Psi(y))<\frac{b_c+m_c^\infty}{2}.
\end{equation}

Now let us define
$$
Q:=\Psi(\overline{B_{\bar{\rho}}(0)}),\ \ \ \ S:=\{u \in H^{s}(\mathbb{R}^{N}): u \in S_{c}\cap P,\ \beta(u)=0\},
$$
$$
\mathcal{H}=\{h \in C(S_{c}\cap P, S_{c}\cap P): h(u)=u,\ \mbox{if}\ I(u)<\frac{b_c+m_c^\infty}{2}\},
$$
we first claim $h(Q) \cap S \neq \emptyset$ for any $h \in \mathcal{H}$. Indeed, given $h \in \mathcal{H}$, we define $T: \overline{B_{\bar{\rho}}(0)} \rightarrow \mathbb{R}^{N}$ for $T(y)=\beta \circ h \circ \Psi(y)$. Because $T$ is the composition of continuous functions, then it is continuous.
Note that for any $|y|=\bar{\rho}$, we have $\Psi(y) \in \partial Q$. Hence from \eqref{IP} , we obtain $h\circ \Psi(y)=\Psi(y)$ and then $T(y)=y$. By the Fixed Point Theorem of Brouwer, we obtain that there exists $\tilde{y} \in B_{\bar{\rho}}(0)$ such that $T(\tilde{y})=0$, which implies $h(\Psi(\tilde{y})) \in S$. Hence, $h(Q) \cap S \neq \emptyset$.

Now define
$$
d=\inf_{h \in \mathcal{H}} \max_{u \in Q} I(h(u)),
$$
we have $d \in(m^{\infty}_{c}, \eta m^{\infty}_{c})$. From \eqref{IP} and $\frac{b_c+m_c^\infty}{2}<b_c$, we have
\begin{equation}\label{eq216}
b_{c}=\inf _{S} I>\max _{\partial Q}I>m^{\infty}_{c}.
\end{equation}
Since $h(Q)\cap S \neq\emptyset$ for any $h\in \mathcal{H}$,
let $h\in \mathcal{H}$ be fixed, there exists $\omega\in S$ such that $\omega$ also belongs to $h(Q)$, which means that $\omega=h(v)$ for some $v \in \Psi(\overline{B_{\bar{\rho}}(0)})$. Hence, we have
$$
I(\omega)\geqslant\inf_{u\in S}I(u)\ \ \text { and }\ \  \max_{u \in Q} I(h(u)) \geqslant I(h(v)) .
$$
This implies
$$
\max_{u\in Q}I(h(u))\geqslant I(h(v))=I(\omega) \geqslant \inf _{u \in S} I(u)=b_{c},
$$
and then
$$
d=\inf_{h\in \mathcal{H}}\max_{u\in Q}I(h(u))\geqslant b_{c},
$$
it follows from \eqref{eq216} that $d>m^{\infty}_{c}$. In particular, let $h=id|_{\mathbb{R}^{N}}$, according to Lemma 3.3, we have
$$
\inf_{h \in \mathcal{H}}\max_{u\in Q} I(h(u))\leq\max _{u \in Q} I(u)<\eta m^{\infty}_{c},
$$
which implies that $d<\eta m^{\infty}_{c}$. Hence, we obtain $d \in(m^{\infty}_{c}, \eta m^{\infty}_{c})$.

At last, we claim for any $\mu\in (0, d-\frac{b_c+m_c^\infty}{2})$, there is $u_\mu\in I^{-1}(d-\mu, d+\mu)\cap S_{c}\cap P$ such that
$$
\|I'(u_\mu)\|_*<\mu,
$$
where
$$
\|I'(u)\|_*=\sup\{|I'(u)v|: v\in T_uS_c \mbox{ and } \|v\|\leq 1\},\ \ \forall\ u\in S_c.
$$

Arguing by contradiction, suppose that there is $\mu_0>0$ with $\mu_0<d-\frac{b_c+m_c^\infty}{2}$ such that
$$
\|I'(u)\|_*\geq\frac{\mu_0}{2},\ \ \forall u\in I^{-1}(d-\mu_0, d+\mu_0)\cap S_{c}\cap P.
$$
Applying a deformation lemma (see Lemma 5.15 in \cite{Willem}), there exists $\tau\in C([0,1]\times S_{c}\cap P,S_{c}\cap P)$ such that
$$
\tau(t,u)=u,\ \mbox{ if } u\not\in I^{-1}(d-\mu_0, d+\mu_0)
$$
and
\begin{equation}\label{DSF}
\tau(1,I^{d+\frac{\mu_0}{2}})\subset I^{d-\frac{\mu_0}{2}},
\end{equation}
where $I^\gamma=\{u\in S_{c}\cap P: I(u)\leq\gamma\}$.

From the definition of $d$, there is $h\in \mathcal{H}$ such that
$$
\max_{u\in Q}I(h(u))\leq d+\frac{\mu_0}{2},
$$
that is, $h(Q)\subset I^{d+\frac{\mu_0}{2}}$. Hence, from \eqref{DSF}, we have
$$
\tau(1,h(Q))\subset I^{d-\frac{\mu_0}{2}},
$$
which implies that
\begin{equation}\label{P}
\max_{u\in \tau(1,h(Q))}I(u)\leq d-\frac{\mu_0}{2}.
\end{equation}

On the other hand, we define
$$
\tilde{h}:=\tau(1,\cdot)\circ h\in C(S_{c}\cap P, S_{c}\cap P).
$$
Note that if $I(u)<\frac{b_c+m_c^\infty}{2}<d-\mu_0$, we have
$$
h(u)=u\ \ \ \mbox{and}\ \ \ \tau(1,u)=u,
$$
thus, $\tilde{h}\in \mathcal{H}$ and $\tau(1,h(Q))=\tilde{h}(Q)$. By the definition of $d$, we have
$$
d\leq \max_{u\in \tau(1,h(Q))}I(u),
$$
which is a contradiction with \eqref{P}.

Now there is $u_n\in I^{-1}(d-\mu_n, d+\mu_n)\cap S_{c}\cap P$ with $\mu_n\to 0$ such that
$$
\|I'(u_n)\|_*<\mu_n,\ \ \forall n\in \mathbb{N}
$$
that is, $\{u_{n}\}$ is a non-negative Palais-Smale sequence for $I|_{S_{c}}$ at level $d$. It follows from Lemma $3.2$ that $\{u_{n}\}$ has a convergent subsequence in  $H^{s}(\mathbb{R}^{N})$. Hence problem \eqref{eq1} has a positive solution $u \in S_{c}$ for some $\lambda<0$.\ \ $\Box$\\

\end{document}